\documentclass[a4paper, 11pt, final]{article}

\usepackage{amsfonts,enumerate}
\usepackage{graphicx,psfrag}
\usepackage{lscape}
\newtheorem{proposition}{Proposition}[section]
\newtheorem{theorem}[proposition]{Theorem}

\newtheorem{definition}[proposition]{Definition}

{\hspace{-5pt}{\nobreak\quad\nobreak\hfill\nobreak$\square$\vspace{8pt}%
\par}\smallskip\goodbreak}

\newenvironment{proofof}[1]{\smallskip\noindent\emph{\textbf{Proof of #1.}}%
\hspace{1pt}}{\hspace{-5pt}{\nobreak\quad\nobreak\hfill\nobreak%
$\square$\vspace{8pt}\par}\smallskip\goodbreak}

\newcommand{\Section}[1]{\section{#1}\setcounter{equation}{0}}

\newcommand{\pint}[1]{\mathaccent23{#1}}
\newcommand{\Lloc}[1]{\mathbf{L^{#1}_{loc}}}
\newcommand{\C}[1]{\mathbf{C^{#1}}}

\newcommand{\PC}{\mathbf{PC}}

\newcommand{\modulo}[1]{{\left|#1\right|}}
\newcommand{\norma}[1]{{\left\|#1\right\|}}

\newcommand{\reali}{{\mathbb{R}}}

\newcommand{\naturali}{{\mathbb{N}}}
\newcommand{\tv}{\mathrm{TV}}
\newcommand{\BV}{\mathbf{BV}}

\renewcommand{\epsilon}{\varepsilon}
\renewcommand{\phi}{\varphi}
\renewcommand{\theta}{\vartheta}
\renewcommand{\L}[1]{\mathbf{L^#1}}
\newcommand{\W}[1]{\mathbf{W^{#1}}}

\begin{document}

\title{\LARGE{Coupling Conditions for the $3\times 3$ Euler System}}

\author{Rinaldo M.~Colombo$^1$ \quad Francesca Marcellini$^2$}
\footnotetext[1]{Dipartimento di Matematica, Universit\`a degli Studi
  di Brescia, Via Branze 38, 25123 Brescia, Italy,
  \texttt{Rinaldo.Colombo@UniBs.it}} \footnotetext[2]{Dipartimento di
  Matematica e Applicazioni, Universit\`{a} di Milano--Bicocca, Via
  Cozzi 53, 20126 Milano, Italy,
  \texttt{F.Marcellini@Campus.Unimib.it}}
\date{}
\maketitle

\begin{abstract}
  This paper is devoted to the extension to the full $3\times3$ Euler
  system of the basic analytical properties of the equations governing
  a fluid flowing in a duct with varying section. First, we consider
  the Cauchy problem for a pipeline consisting of 2 ducts joined at a
  junction. Then, this result is extended to more complex pipes. A key
  assumption in these theorems is the boundedness of the total
  variation of the pipe's section. We provide explicit examples to
  show that this bound is necessary.


  \smallskip\noindent\textit{2000~MSC:} 35L65, 76N10
\end{abstract}

\section{Introduction}

We consider Euler equations for the evolution of a fluid flowing in a
pipe with varying section $a = a(x)$, see~\cite[Section~8.1]{Whitham}
or~\cite{GoatinLefloch2004, LiuDuct2}:
\begin{equation}
  \label{1eq:E1}
  \left\{ 
    \begin{array}{ll}
      \displaystyle
      \partial_{t} (a\rho) 
      + 
      \partial_{x} (a q) 
      = 
      0 
      \\ 
      \displaystyle
      \partial_{t} (aq) 
      +
      \partial_{x} 
      \left[ 
        aP\left( \rho,q,E \right)
      \right]  
      =
      p\left( \rho,e\right)  \, \partial_{x}a
      \\ 
      \displaystyle
      \partial_{t} (aE) 
      + 
      \partial_{x} \left[ a F\left( \rho,q,E\right) \right]  
      = 
      0
    \end{array}
  \right.  
\end{equation}
where, as usual, $\rho$ is the fluid density, $q$ is the linear
momentum density and $E$ is the total energy density. Moreover
\begin{equation}
  \label{1eq:EPF}
 \!\!\!
 E(\rho, q, E) = \frac{1}{2}\frac{q^{2}}{\rho}+\rho e
  ,\quad 
  P(\rho, q, E) = \frac{q^{2}}{\rho}+p
  , \quad 
  F(\rho, q, E) = \frac{q}{\rho} ( E + p) \,,
\end{equation}
with $e$ being the internal energy density, $P$ the flow of the linear
momentum density and $F$ the flow of the energy density. The above
equations express the conservation laws for the mass, momentum, and
total energy of the fluid through the pipe. Below, we will often refer
to the standard case of the ideal gas, characterized by the relations
\begin{equation}
  \label{1eq:perfect}
  p = (\gamma - 1) \rho e
  , \qquad 
  S = \ln e - ( \gamma - 1) \ln\rho \,,
\end{equation}
for a suitable $\gamma>1$. Note however, that this particular equation
of state is necessary only in case~\textbf{(p)} of
Proposition~\ref{1prop:otherpsi} and has been used in the examples in
Section~\ref{1sec:Blow-Up}. In the rest of this work, the usual
hypothesis~\cite[formula~(18.8)]{Smoller}, that is $p>0$,
$\partial_\tau p(\tau,S) < 0$ and $\partial^2_{\tau\tau} p(\tau,S)>0$,
are sufficient.

The case of a sharp discontinuous change in the pipe's section due to
a junction sited at, say, $x=0$, corresponds to $a(x) = a^-$ for $x <
0$ and $a(x) = a^+$ for $x > 0$. Then, the motion of the fluid can be
described by
\begin{equation}
  \label{1eq:E}
  \left\{ 
    \begin{array}{ll}
      \displaystyle
      \partial_{t} \rho 
      + 
      \partial_{x} q 
      = 
      0
      \\ 
      \displaystyle
      \partial_{t} q 
      +
      \partial_{x} 
      P (\rho,q,E )
      =
      0
      \\ 
      \displaystyle
      \partial_{t}E 
      + 
      \partial_{x}  F(\rho,q,E)
      = 
      0,
    \end{array}
  \right.
\end{equation}
for $x \neq 0$, together with a \emph{coupling condition} at the
junction of the form:
\begin{equation}
  \label{1eq:Junction}
  \Psi \left( a^-,(\rho,q,E)(t, 0-); a^+, (\rho,q,E)(t, 0+) \right) = 0.
\end{equation}
Above, we require the existence of the traces at $x=0$ of $(\rho, q,
E)$. Various choices of the function $\Psi$ are present in the
literature, see for instance~\cite{BandaHertyKlar1,
  ColomboGaravello2009, ColomboHertySachers,
  ColomboMarcellini} in the case of the $p$-system
and~\cite{ColomboMauri} for the full $3\times 3$
system~(\ref{1eq:E}). Here, we consider the case of a general coupling
condition which comprises all the cases found in the
literature. Within this setting, we prove the well posedness of the
Cauchy problem for~(\ref{1eq:E})--(\ref{1eq:Junction}). Once this result
is obtained, the extension to pipes with several junctions and to
pipes with a $\W{1,1}$ section is achieved by the standard methods
considered in the literature. For the analytical techniques to cope
with networks having more complex geometry, we refer
to~\cite{GaravelloPiccoliBook}.

The above statements are global in time and local in the space of the
thermodynamic variables $(\rho, q, E)$. Indeed, for any fixed
(subsonic) state $(\bar \rho, \bar q, \bar E)$, there exists a bound
on the total variation $\tv(a)$ of the pipe's section, such that all
sections below this bound give rise to Cauchy problems
for~(\ref{1eq:E})--(\ref{1eq:Junction}) that are well posed in $\L1$. We
show the necessity of this bound in the conditions found in the
current literature. Indeed, we provide explicit examples showing that
a wave can be arbitrarily amplified through consecutive interactions
with the pipe walls, see Figure~\ref{1fig:f}.

The paper is organized as follows. The next section is divided into
three parts, the former one deals with a single junction and two
pipes, then we consider $n$ junctions and $n+1$ pipes, the latter part
presents the case of a $\W{1,1}$ section. Section~\ref{1sec:Coupl Cond}
is devoted to different specific choices of coupling
conditions~(\ref{1eq:Junction}). In Section~\ref{1sec:Blow-Up}, an
explicit example shows the necessity of the bound on the total
variation of the pipe's section. All proofs are gathered in
Section~\ref{1sec:Tech}.

\section{Basic Well Posedness Results}
\label{1sec:General}

Throughout, we let $u = (\rho, q, E)$. We denote by $\reali^{+}$ the
real halfline $\left[ 0,+\infty \right[$, while $\pint{\reali}^{+} =
\left] 0, +\infty \right[ $.  Following various results in the
literature, such as~\cite{BandaHertyKlar1, BandaHertyKlar2,
  ColomboGaravello2009, ColomboHertySachers, ColomboMarcellini,
  ColomboMauri, GuerraMarcelliniSchleper}, we limit the analysis in
this paper to the \textit{subsonic} region given by $\lambda_1(u) < 0
< \lambda_3(u)$ and $\lambda_2(u) \neq 0$, where $\lambda_{i}$ is the
$i-$th eigenvalue of~(\ref{1eq:E}), see~(\ref{1eigenv}). Without any
loss of generality, we further restrict to
\begin{equation}
  \label{1subsonic}
  A_{0}
  = 
  \left\{ 
    u \in \pint\reali^{+} \times \reali^{+}\times \pint\reali^{+}
    \colon 
    \lambda_{1} (u) < 0 < \lambda_{2} (u)
  \right\} \,.
\end{equation}
Note that we fix \emph{a priori} the sign of the fluid speed $v$,
since $\lambda_{2} (u) = q / \rho = v > 0$.

\subsection{A Junction and two Pipes}

We now give the definition of weak $\Psi-$solution to the Cauchy
Problem for~(\ref{1eq:E}) equipped with the
condition~(\ref{1eq:Junction}),
extending~\cite[Definition~2.1]{ColomboGaravello2009}
and~\cite[Definition~2.2]{ColomboMarcellini} to the $3\times 3$
case~(\ref{1eq:E}) and comprising the particular case covered
in~\cite[Definition~2.4]{ColomboMauri}.

\begin{definition}
  \label{1def:Weaksolution}
  Let $\Psi \colon (\pint{\reali}^+ \times A_0)^2 \to \reali^3$, $u_o
  \in \BV(\reali; A_0)$ and two positive sections $a^-$, $a^+$ be
  given.  A $\Psi$-solution to~(\ref{1eq:E}) with initial datum $u_o$
  is a map
  \begin{equation}
    \label{1eq:regularity}
    \begin{array}{rcl}
      u 
      & \in & 
      \C0 \left( \reali^+; \Lloc1 (\reali^{+}; A_0 ) \right)
      \\ 
      u (t) & \in & \BV (\reali; A_0 )
      \quad \mbox{ for a.e. } t \in \reali^+
    \end{array}
  \end{equation}
  such that
  \begin{description}
  \item[1.] for $x \neq 0$, $u$ is a weak entropy solution
    to~(\ref{1eq:E});
  \item[2.] for a.e.~$x\in \reali$, $u(0,x) = u_o(x)$;
  \item[3.] for a.e.~$t\in \reali^+$, the coupling condition
    (\ref{1eq:Junction}) at the junction is met.
  \end{description}
\end{definition}

\noindent Below, extending the $2\times2$ case of the $p$-system,
see~\cite{BandaHertyKlar1, ColomboGaravello2007, ColomboGaravello2009,
  ColomboHertySachers, ColomboMarcellini}, we consider some properties
of the coupling condition (\ref{1eq:Junction}), which we rewrite here
as
\begin{equation}
  \label{1eq:Psi}
  \Psi(a^-,u^-;a^+,u^+)
  =0\,.
\end{equation}
\begin{description}
\item[($\mathbf{\Psi}$0)] Regularity: $ \Psi \in \C1 \left( (
    \pint{\reali}^+ \times A_0) ^{2} ; \reali^3 \right)$.
\item[($\mathbf{\Psi}$1)] No-junction case: for all $a >0$ and all
  $u^-,u^{+} \in A_0 $, then
  \begin{displaymath}
    \Psi(a,u^-;a,u^+) = 0  \Longleftrightarrow u^- = u^{+} \,.
  \end{displaymath}

\item[($\mathbf{\Psi}$2)] Consistency: for all positive $a^-, a^0,a^+$
  and all $u^-,u^{0},u^{+}\in A_0$,
  \begin{displaymath}
    \begin{array}{l}
      \Psi(a^{-},u^-;a^{0},u^0) = 0
      \\
      \Psi(a^{0},u^{0};a^{+},u^+) = 0
    \end{array}
    \Longrightarrow
    \Psi(a^{-},u^-;a^{+},u^+)=0 \,.
  \end{displaymath}
\end{description}

\noindent Moreover, by an immediate extension
of~\cite[Lemma~2.1]{ColomboMarcellini}, ($\mathbf{\Psi0}$) ensures
that (\ref{1eq:Psi}) implicitly defines a map
\begin{equation}
  \label{1eq:T}
  u^{+} = T ( a^{-} , a^{+} ;u^{-}) 
\end{equation}
in a neighborhood of any pair of subsonic states $u^-$, $u^+$ and
sections $a^-, a^+$ that satisfy $\Psi(a^-,u^-;a^+,u^+) = 0$.

The technique in~\cite{ColomboGuerra6} allows to prove the following
well posedness result.

\begin{theorem}
  \label{1thm:CauchyProblem}
  Assume that $\Psi$ satisfies
  conditions~(\textbf{$\mathbf{\Psi}$0})-(\textbf{$\mathbf{\Psi}$2}). For
  every $\bar a > 0$ and $\bar u \in A_0$ such that
  \begin{equation}
    \label{1eq:det}
    \det \left[ 
      D_{u^{-}}\Psi\cdot r_{1}(\bar u)
      \quad
      D_{u^{+}}\Psi \cdot r_{2}(\bar u)
      \quad
      D_{u^{+}}\Psi \cdot r_{3}(\bar u)
    \right]\neq 0
  \end{equation}
  there exist positive $\delta$, $L$ such that for all $a^-, a^+$ with
  $\modulo{a^+ - \bar a} + \modulo{a^- - \bar a} < \delta$ there
  exists a semigroup $S \colon \reali^ + \times \mathcal{D} \to
  \mathcal{D}$ with the following properties:
  \begin{enumerate}
  \item $\mathcal{D} \supseteq \left\{ u \in \bar u + \L1(\reali; A_0)
      \colon \tv(u) < \delta \right\}$.
  \item For all $u \in \mathcal{D}$, $S_0 u = u$ and for all $t,s \geq
    0$, $S_t S_s u = S_{s+t} u$.
  \item For all $u, u' \in \mathcal{D}$ and for all $t,t' \geq 0$,
    \begin{displaymath}
      \norma{S_t u - S_{t'} u'}_{\L1}
      \leq
      L \cdot \left(
        \norma{u - u'}_{\L1} + \modulo{t-t'}
      \right)
    \end{displaymath}
  \item If $u \in \mathcal{D}$ is piecewise constant, then for $t$
    small, $S_t u$ is the gluing of solutions to Riemann problems at
    the points of jump in $u$ and at the junction at $x = 0$.
  \item For all $u_o \in \mathcal{D}$, the orbit $t \to S_t u_o$ is a
    $\Psi$-solution to~(\ref{1eq:E}) with initial datum $u_o$.
  \end{enumerate}
\end{theorem}

\noindent The proof is postponed to Section~\ref{1sec:Tech}. Above
$r_{i}( u)$, with $i=1,2,3$, are the right eigenvectors of $Df(u)$,
see~(\ref{1eigenv}). Moreover, by \emph{solution to the Riemann
  Problems at the points of jump} we mean the usual Lax solution,
see~\cite[Chapter~5]{BressanLectureNotes}, whereas for the definition
of \emph{solution to the Riemann Problems at the junction} we refer
to~\cite[Definition~2.1]{ColomboHertySachers}.

\subsection{$n$ Junctions and $n+1$ Pipes}

The same procedure used in~\cite[Paragraph~2.2]{ColomboMarcellini}
allows now to construct the semigroup generated by~(\ref{1eq:E}) in the
case of a pipe with piecewise constant section
\begin{displaymath}
  a 
  = 
  a_0 \, \chi_{]-\infty, x_1] }
  + 
  \sum_{j=1}^{n-1} a_j \, \chi_{[x_{j}, x_{j+1}[}
  +
  a_{n} \, \chi_{[x_n, +\infty[}
\end{displaymath}
with $n \in \naturali$. In each segment $\left]x_j, x_{j+1} \right[$,
the fluid is modeled by~(\ref{1eq:E}). At each junction $x_j$, we
require condition~(\ref{1eq:Junction}), namely
\begin{equation}
  \label{1eq:Psii}
  \Psi (a_{j-1}, u_j^-; a_{j}, u_j^+) 
  =
  0
  \quad 
  \begin{array}{l}
    \mbox{ for all } j = 1, \ldots, n \mbox{, where}
    \\
    \displaystyle  u_j^\pm = \lim_{x \to x_j\pm} u_j (x) \,.
  \end{array}
\end{equation}
We omit the formal definition of $\Psi$-solution
to~(\ref{1eq:E})--(\ref{1eq:Junction}) in the present case, since it is
an obvious iteration of Definition~\ref{1def:Weaksolution}. The natural
extension of Theorem~\ref{1thm:CauchyProblem} to the case
of~(\ref{1eq:E})--(\ref{1eq:Psii}) is the following result.

\begin{theorem}
  \label{1thm:nn}
  Assume that $\Psi$ satisfies
  conditions~(\textbf{$\mathbf{\Psi}$0})-(\textbf{$\mathbf{\Psi}$2}). For
  any $\bar a>0$ and any $\bar u \in A_0$, there exist positive $M,
  \Delta,\delta, L, \mathcal{M}$ such that for any pipe's profile
  satisfying
  \begin{equation}
    \label{1eq:A0}
    a \in \PC \left(\reali; \left]\bar a- \Delta, \bar
        a+\Delta\right[\right) \mbox{ with } \tv(a) < M
  \end{equation}
  \noindent there exists a piecewise constant stationary solution
  \begin{displaymath}
    \hat u
    = 
    \hat u_0 \chi_{\left]-\infty, x_1\right[}
    +
    \sum_{j=1}^{n-1} \hat u_j \chi_{\left]x_{j},  x_{j+1}\right[}
    +
    \hat u_n \chi_{\left]x_{n}, +\infty\right[}
  \end{displaymath}
  to~(\ref{1eq:E})--(\ref{1eq:Psii}) satisfying
  \begin{eqnarray}
    \nonumber
    & &
    \hat u_j \in A_0 \mbox{ with }
    \modulo{\hat u_j - \bar u} < \delta
    \mbox{ for }
    j=0, \ldots n
    \\
    \nonumber
    & &
    \Psi
    \left( 
      a_{j-1}, \hat u_{j-1}; a_{j}, \hat u_{j}
    \right)
    =
    0
    \mbox{ for }
    j = 1, \ldots, n
    \\
    \label{1eq:tvhat}
    & &
    \tv(\hat u) \leq \mathcal{M} \, \tv(a)
  \end{eqnarray}
  and a semigroup $S^a \colon \reali^ + \times \mathcal{D}^a \to
  \mathcal{D}^a$ such that
  \begin{enumerate}
  \item $\mathcal{D}^a \supseteq \left\{ u \in \hat u + \L1(\reali;
      A_0) \colon \tv(u - \hat u) < \delta \right\}$.
  \item $S^{a}_0$ is the identity and for all $t,s \geq 0$, $S^a_t
    S^a_s = S^a_{s+t}$.
  \item For all $u, u' \in \mathcal{D}^a$ and for all $t, t' \geq 0$,
    \begin{displaymath}
      \norma{S^a_t u - S^a_{t'} u'}_{\L1}
      \leq
      L \cdot \left(
        \norma{(u) - u'}_{\L1} + \modulo{t-t'}
      \right) .
    \end{displaymath}
  \item If $u \in \mathcal{D}^a$ is piecewise constant, then for $t$
    small, $S_t u$ is the gluing of solutions to Riemann problems at
    the points of jump in $u$ and at each junction $x_j$.
  \item For all $u \in \mathcal{D}^a$, the orbit $t \to S^a_t u$ is a
    weak $\Psi$-solution to~(\ref{1eq:E})--(\ref{1eq:Psii}).
  \end{enumerate}
\end{theorem}

\noindent We omit the proof, since it is based on the natural
extension to the present $3\times 3$ case
of~\cite[Theorem~2.4]{ColomboMarcellini}.  Remark that, as in that
case, $\delta$ and $L$ depend on $a$ only through $\bar a$ and
$\tv(a)$. In particular, all the construction above is independent
from the number of points of jump in $a$.

\subsection{A Pipe with a $\W{1,1}$ Section}

In this paragraph, the pipe's section $a$ is assumed to satisfy
\begin{equation}
  \label{1eq:A1}
  \left\{
    \begin{array}{l}
      a \in \W{1,1} \left(\reali;\left]\bar a-\Delta, \bar a+\Delta\right[
      \right) \mbox{ for suitable } \Delta >0, \ \bar a > \Delta
      \\
      \tv(a) < M \mbox{ for a suitable } M > 0
      \\
      a'(x) = 0 \mbox{ for a.e. } 
      x \in \reali \setminus [-X, X]
      \mbox{ for a suitable } X > 0 \,.
    \end{array}
  \right.
  \!\!
\end{equation}
The same procedure used in~\cite[Theorem~2.8]{ColomboMarcellini}
allows to construct the semigroup generated by~(\ref{1eq:E1}) in the
case of a pipe which satisfies~(\ref{1eq:A1}). Indeed, thanks to
Theorem~\ref{1thm:nn}, we approximate $a$ with a piecewise constant
function $a_n$. The corresponding
problems~to~(\ref{1eq:E})--(\ref{1eq:Psii}) generate semigroups $S_n$
defined on domains characterized by uniform bounds on the total
variation and that are uniformly Lipschitz in time. Here, uniform
means also independent from the number of junctions. Therefore, we
prove the pointwise convergence of the $S_n$ to a limit semigroup $S$,
along the same lines in~\cite[Theorem~2.8]{ColomboMarcellini}.

\Section{Coupling Conditions}
\label{1sec:Coupl Cond}
This section is devoted to different specific choices
of~(\ref{1eq:Psi}).

\paragraph{\textbf{(S)-Solutions}}

We consider first the coupling condition inherited from the smooth
case.  For smooth solutions and pipes' sections, system~(\ref{1eq:E1})
is equivalent to the $3\times 3$ balance law
\begin{equation}
  \label{1eq:E2}
  \left\{ 
    \begin{array}{ll}
      \displaystyle
      \partial_{t} \rho 
      + 
      \partial_{x}q 
      = 
      -\frac{q}{a} \, \partial _{x} a 
      \\ 
      \displaystyle
      \partial_{t} q 
      +
      \partial_{x} 
      P(\rho,q,E)
      =
      -\frac{q^{2}}{a\rho } \, \partial_{x}a
      \\ 
      \displaystyle
      \partial_{t}E 
      + 
      \partial_{x}  F(\rho,q,E)
      = 
      -\frac{F}{a} \, \partial _{x} a.
    \end{array}
  \right.
\end{equation}
The stationary solutions to~(\ref{1eq:E1}) are characterized as
solutions to
\begin{equation}
  \label{1eq:Stationary}
  \!\!\!
  \left\{
    \!\!\!\!
    \begin{array}{ll}
      \partial_x (a(x) \, q) = 0
      \\
      \partial_x \!
      \left( 
        a(x) \,P(\rho,q,E)
      \right) 
      =
      p(\rho, e)\, \partial_x a
      \\
      \partial_x \!
      \left( 
        a(x) \,F(\rho,q,E)
      \right) 
      =0
    \end{array}
  \right.
  \mbox{ or }
  \left\{ 
    \!\!\!\!
    \begin{array}{ll}
      \partial_{x} q = 
      -\frac{q}{a} \, \partial _{x} a 
      \\ 
      \partial_{x} 
      P\left( \rho,q,E \right) =
      -\frac{q^{2}}{a\rho } \, \partial_{x}a
      \\ 
      \partial_{x} 
      F\left( \rho,q,E \right) =
      -\frac{F}{a } \, \partial_{x}a \,.
    \end{array}
  \right.  
  \!\!\!
\end{equation}
As in the $2\times2$ case of the $p$-system, the smoothness of the
sections induces a unique choice for condition~(\ref{1eq:Psi}),
see~\cite[(2.3) and~(2.19)]{ColomboMarcellini}, which reads
\begin{equation}
  \label{1eq:psismooth}
  \mbox{\textbf{(S)}}  
  \quad
  \Psi 
  =
  \left[\!\!\!
    \begin{array}{l}
      a^+ q^+ - a^- q^-
      \\
      \displaystyle
      a^+ P(u^+) - a^- P(u^-) 
      +
      \int_{-X}^X p \left( \mathcal{R}^{a}(x), \mathcal{E}^a(x) \right) 
      a'(x)
      \mathrm{d}x\!
      \\
      a^+ F(u^+) - a^- F(u^-) 
    \end{array}\!
  \right]
\end{equation}
where $a=a(x)$ is a smooth monotone function satisfying $a(-X) = a^-$
and $a(X)=a^+$, for a suitable $X>0$. $\mathcal{R}^a,\mathcal{E}^a$
are the $\rho$ and $e$ component in the solution
to~(\ref{1eq:Stationary}) with initial datum $u^{-}$ assigned at
$-X$. Note that, by the particular form of~(\ref{1eq:psismooth}), the
function $\Psi$ is independent both from the choice of $X$ and from
that of the map $a$, see~\cite[2.~in
Proposition~2.7]{ColomboMarcellini}.

\paragraph{\textbf{(P)-Solutions}}

The particular choice of the coupling condition
in~\cite[Section~3]{ColomboMauri} can be recovered in the present
setting. Indeed, conditions~\textbf{(M)}, \textbf{(E)}
and~\textbf{(P)} therein amount to the choice
\begin{equation}
  \label{1eq:psiCri}
  \mbox{\textbf{(P)}}  
  \qquad
  \Psi(a^{-},u^{-},a^{+},u^{+})
  =
  \left[
    \begin{array}{c}
      a^{+}  q^{+}-a^{-} q^{-}
      \\
      P(u^+) -  P(u^-)
      \\
      a^+ F(u^+) - a^- F(u^-)
    \end{array}
  \right],
\end{equation}
where $a^{+}$ and $a^{-}$ are the pipe's sections. Consider fluid
flowing in a horizontal pipe with an elbow or kink,
see~\cite{HoldenRisebroKink}. Then, it is natural to assume the
conservation of the total linear momentum along directions dependent
upon the geometry of the elbow. As the angle of the elbow vanishes,
one obtains the condition above,
see~\cite[Proposition~2.6]{ColomboMauri}.

\paragraph{\textbf{(L)-Solutions}}

We can extend the construction in~\cite{BandaHertyKlar1,
  BandaHertyKlar2, ColomboGaravello2007} to the $3\times 3$
case~(\ref{1eq:E}). Indeed, the conservation of the mass and linear
momentum in~\cite{ColomboGaravello2007} with the conservation of the
total energy for the third component lead to the choice
\begin{equation}
  \label{1eq:psiL}
  \mbox{\textbf{(L)}}  
  \qquad
  \Psi(a^{-},u^{-},a^{+},u^{+})
  =
  \left[
    \begin{array}{c}
      a^{+} q^{+}-a^{-} q^{-}
      \\
      a^+ P(u^+) - a^- P(u^-)
      \\
      a^+ F(u^+) - a^- F(u^-)
    \end{array}
  \right],
\end{equation}
where $a^{+}$ and $a^{-}$ are the pipe's sections. The above is the
most immediate extension of the standard definition of Lax solution to
the case of the Riemann problem at a junction.

\paragraph{\textbf{(p)-Solutions}}

Following~\cite{BandaHertyKlar1,BandaHertyKlar2}, motivated by the
what happens at the hydrostatic equilibrium, we consider a coupling
condition with the conservation of the pressure $p(\rho)$ in the
second component of $\Psi$. Thus

\begin{equation}
  \label{1eq:psip}
  \mbox{\textbf{(p)}}  
  \qquad
  \Psi(a^{-},u^{-},a^{+},u^{+})
  =
  \left[
    \begin{array}{c}
      a^{+} \, q^{+}-a^{-} \, q^{-}
      \\
      p(\rho^{+}, e^{+})-p(\rho^{-}, e^{-})
      \\
      a^+ F(u^+) - a^- F(u^-)
    \end{array}
  \right],
\end{equation}
where $a^{+}$ and $a^{-}$ are the pipe's sections.

\smallskip

\begin{proposition}
  \label{1prop:otherpsi}
  For every $\bar a>0$ and $\bar u\in A_{0}$, each of the coupling
  conditions $\Psi$ in~(\ref{1eq:psismooth}), (\ref{1eq:psiCri}),
  (\ref{1eq:psiL}), (\ref{1eq:psip}) satisfies the
  requirements~(\textbf{$\mathbf{\Psi}$0})-(\textbf{$\mathbf{\Psi}$2})
  and~(\ref{1eq:det}). In the case of~(\ref{1eq:psip}), we also require
  that the fluid is perfect, i.e.~that~(\ref{1eq:perfect}) holds.
\end{proposition}

The proof is postponed to Section~\ref{1sec:Tech}. Thus,
Theorem~\ref{1thm:CauchyProblem} applies, yielding the well posedness
of~(\ref{1eq:E})--(\ref{1eq:Junction}) with each of the particular
choices of $\Psi$ in~(\ref{1eq:psismooth}), (\ref{1eq:psiCri}),
(\ref{1eq:psiL}), (\ref{1eq:psip}).

\Section{Blow-Up of the Total Variation}
\label{1sec:Blow-Up}

In the previous results a key role is played by the bound on the total
variation $\tv(a)$ of the pipe's section. This requirement is
intrinsic to problem~(\ref{1eq:E})--(\ref{1eq:Junction}) and not due to
the technique adopted above. Indeed, we show below that in each of the
cases~(\ref{1eq:psismooth}), (\ref{1eq:psiCri}), (\ref{1eq:psiL}),
(\ref{1eq:psip}), it is possible to choose an initial datum and a
section $a \in \BV(\reali; [a^-, a^+])$ with $a^+ - a^-$ arbitrarily
small, such that the total variation of the corresponding solution
to~(\ref{1eq:E})--(\ref{1eq:Junction}) becomes arbitrarily large.

\begin{figure}[htpb]
  \centering
  \begin{psfrags}
    \psfrag{a}{$x$} \psfrag{e}{$\Delta a$}
    \psfrag{b}{$a$}\psfrag{c}{$x$} \psfrag{f}{$\sigma_3^-$}
    \psfrag{g}{$\sigma_3^+$} \psfrag{h}{$\sigma_3^{++}$}
    \psfrag{d}{$t$} \psfrag{l}{$u^+$} \psfrag{i}{$u$} \psfrag{n}{$2l$}
    \psfrag{m}{$l$}
    \includegraphics[width=6cm]{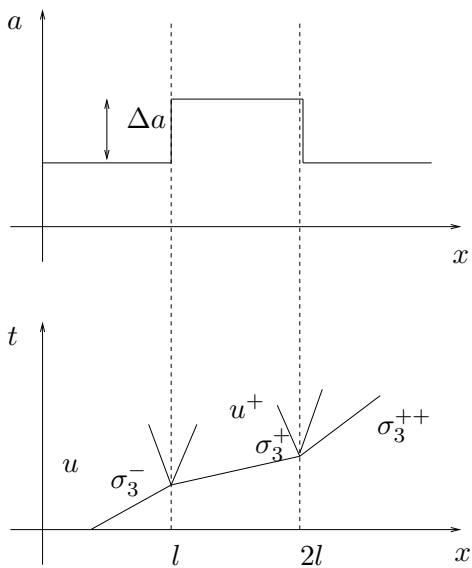}
  \end{psfrags}
  \caption{\label{1fig:f}A wave $\sigma_3^-$ hits a junction where the
    pipe's section increases by $\Delta a$. From this interaction, the
    wave $\sigma_3^+$ arises, which hits a second junction, where the
    pipe section decreases by $\Delta a$.}
\end{figure}

Consider the case in Figure~\ref{1fig:f}.
A wave $\sigma_3^-$ hits a junction where the pipe's section increases
by, say, $\Delta a>0$. The fastest wave arising from this interaction
is $\sigma_3^+$, which hits the second junction where the section
diminishes by $\Delta a$.

\begin{figure}[htpb]
  \centering
  \begin{psfrags}
    \psfrag{s3m}{$\!\!\!\!\sigma_3^-$} \psfrag{s1p}{$\sigma_1^+$}
    \psfrag{s2p}{$\sigma_2^+$} \psfrag{s3p}{$\sigma_3^+$}
    \psfrag{1}{$u$} \psfrag{2}{} \psfrag{3}{} \psfrag{4}{$u^+$}
    \psfrag{5}{} \psfrag{6}{}
    \includegraphics[width=4cm]{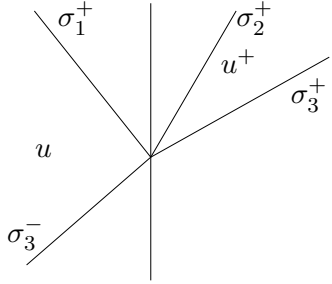}
  \end{psfrags}
  \caption{Notation used in~(\ref{1eq:PbR}) and~(\ref{1second}).}
  \label{1fig:NotInter}
\end{figure}

Solving the Riemann problem at the first interaction amounts to solve
the system
\begin{equation}
  \label{1eq:PbR}
  L_{3}\left( L_2 \left( 
      T \left( 
        L_1 (  u; \sigma_1^+) 
      \right)
      ; \sigma_2^+\right) ;\sigma_3^+
  \right) 
  = 
  T \left(
    L_3( u; \sigma_3^-)
  \right),
\end{equation}
where $u \in A_0$, see~Figure~\ref{1fig:NotInter} for the definitions
of the waves' strengths $\sigma_i^+$ and $\sigma_3^-$. Above, $T$ is
the map defined in~(\ref{1eq:T}), which in turn depends from the
specific condition~(\ref{1eq:Psi}) chosen. In the expansions below, we
use the $(\rho,q,e)$ variables, thus setting $u = (\rho,q,e)$
throughout this section.
Differently from the case of the $2\times2$ $p$-system
in~\cite{ColomboMarcellini}, here we need to consider the second order
expansion in $\Delta a = a^+ - a^-$ of the map $T$; that is
\begin{equation}
  \label{1eq:TT}
  T(a, a+\Delta a; u)
  = 
  u
  +
  H(u) \frac{\Delta a}{a}
  +
  G(u) \left( \frac{\Delta a}{a} \right)^2
  +
  o\left( \frac{\Delta a}{a} \right)^2
\end{equation}
The explicit expressions of $H$ and $G$ in~(\ref{1eq:TT}), for each of
the coupling conditions~(\ref{1eq:psismooth}), (\ref{1eq:psiCri}),
(\ref{1eq:psiL}), (\ref{1eq:psip}), are in Section~\ref{1subs:Conti}.

Inserting~(\ref{1eq:TT}) in the first order expansions in the wave's
sizes of~(\ref{1eq:PbR}), with $\tilde r_{i}$ for $i=1,2,3$ as
in~(\ref{1tangent vectors}), we get a linear system in
$\sigma_{1}^{+},\sigma_{2}^{+},\sigma_{3}^{+}$. Now, introduce the
fluid speed $v = q/\rho$ and the adimensional parameter
\begin{displaymath}
  \theta 
  =
  \left(\frac{v}{c}\right)^2 
  =
  \frac{v^{2}}{\gamma(\gamma-1)e} \,,
\end{displaymath} 
a sort of \emph{``Mach number''}.  Obviously, $\theta \in [0,1]$ for
$u \in A_0$.  We thus obtain an expression for $\sigma_{3}^{+}$ of the
form
\begin{equation}
  \label{1first}
  \sigma_{3}^{+}
  =
  \left( 
    1 
    + 
    f_{1}(\theta)  \, \frac{\Delta a}{a}
    +
    f_{2}(\theta)  \, \left( \frac{\Delta a}{a}\right)^{2}
  \right)
  \sigma_{3}^{-} \,.
\end{equation}
The explicit expressions of $f_1$ and $f_2$ in~(\ref{1first}) are in
Section~\ref{1subs:Conti}.

Remark that the present situation is different from that of the
$2\times2$ $p$-system considered in~\cite{ColomboMarcellini}. Indeed,
for the $p$-system $f_{2}(\theta) = f_{2}(\theta^{+}) = 0$, while here
it is necessary to compute the second order term in $(\Delta a)/a$.

Concerning the second junction, similarly, we introduce the parameter
$\theta^{+} = (v^+/c^+)^2$ which corresponds to the state
$u^{+}$. Recall that $u^+$ is defined by $u^+ = L_{3}^{-} \left( T
  \left( L_3 ( u; \sigma_3^-) ;\sigma_3^+ \right) \right)$, see
Figure~\ref{1fig:NotInter} and Section~\ref{1subs:Conti} for the
explicit expressions of $\theta^+$.  We thus obtain the estimate
\begin{equation}
  \label{1second}
  \sigma_{3}^{++}
  =
  \left( 
    1 
    - 
    f_{1}(\theta^{+})  \, \frac{\Delta a}{a}
    +
    f_{2}(\theta^{+})  \, \left( \frac{\Delta a}{a}\right)^{2}
  \right)
  \sigma_{3}^{+},
\end{equation}
where $\theta^+ = \theta^+ \left(\theta,\sigma_3^-,(\Delta a)/a
\right)$. Now, at the second order in $(\Delta a)/a$ and at the first
order in $\sigma_3^-$, (\ref{1first}) and (\ref{1second}) give
\begin{eqnarray}
  \nonumber
  \sigma_{3}^{++}
  & = &
  \left( 
    1 
    - 
    f_{1}( \theta^{+})  \, \frac{\Delta a}{a}
    +
    f_{2}( \theta^{+})  \, \left( \frac{\Delta a}{a}\right)^{2}
  \right)
  \\
  \nonumber
  & &
  \qquad
  \times
  \left(
    1 
    + 
    f_{1}(\theta)  \, \frac{\Delta a}{a}
    +
    f_{2}(\theta)  \, \left( \frac{\Delta a}{a}\right)^{2}
  \right)
  \sigma_{3}^{-} 
  \\
  & = &
  \label{1eq:chi}
  \left( 1 + \chi(\theta) \left(  \frac{\Delta a}{a}\right) ^{2} \right) 
  \sigma_{3}^{-} \,.
\end{eqnarray}
Indeed, computations show that $f_{1}\left( \theta\right)-f_{1}\left(
  \theta^{+}\right)$ vanishes at the first order in $(\Delta a)/a$, as
in the case of the $p$-system. The explicit expressions of $\chi$ are
in Section~\ref{1subs:Conti}.

It is now sufficient to compute the sign of $\chi$. If it is positive,
then repeating the interaction in Figure~\ref{1fig:f} a sufficient
number of times leads to an arbitrarily high value of the refracted
wave $\sigma_3$ and, hence, of the total variation of the solution
$u$.

Below, Section~\ref{1sec:Tech} is devoted to the computations of $\chi$
in the different cases~(\ref{1eq:psismooth}), (\ref{1eq:psiCri}),
(\ref{1eq:psiL}) and~(\ref{1eq:psip}). To reduce the formal complexities
of the explicit computations below, we consider the standard case of
an ideal gas characterized by~(\ref{1eq:perfect}) with $\gamma = 5 /
3$.

\begin{figure}[htpb]
  \includegraphics[width=0.49\linewidth]{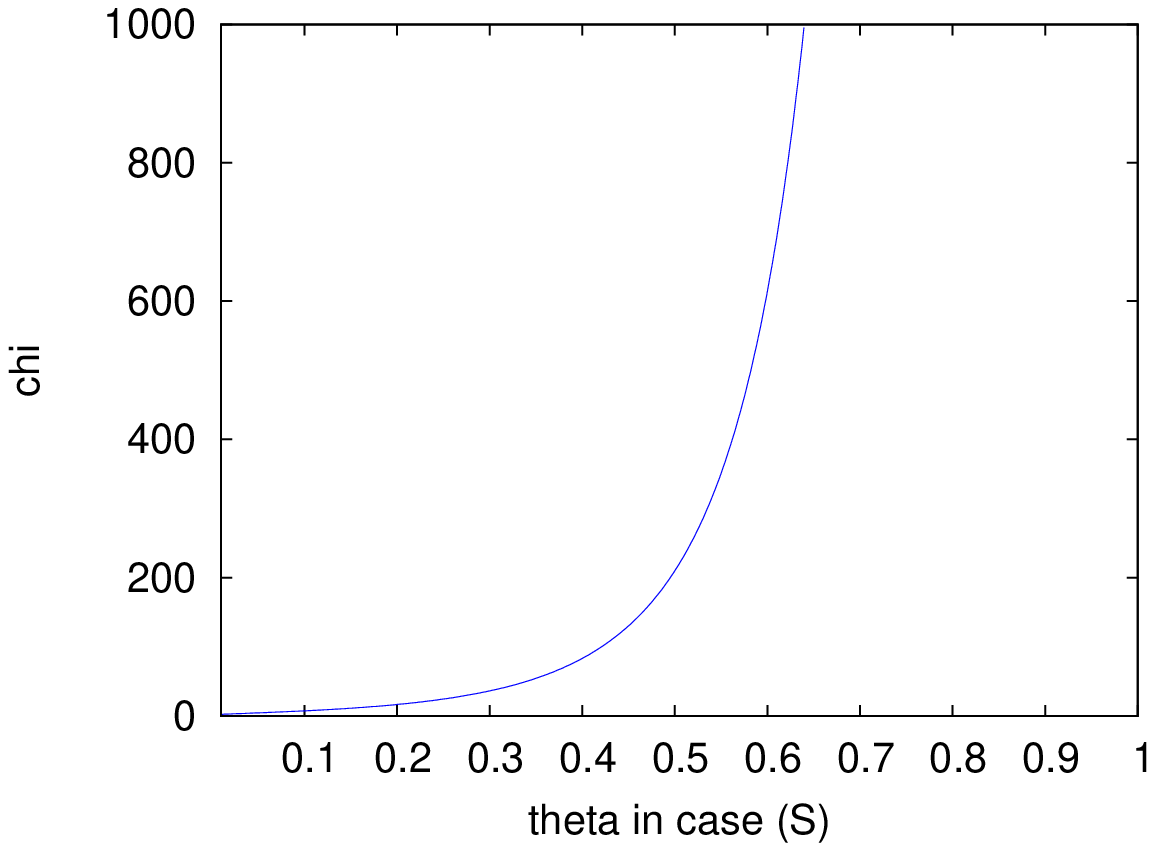}
  \includegraphics[width=0.49\linewidth]{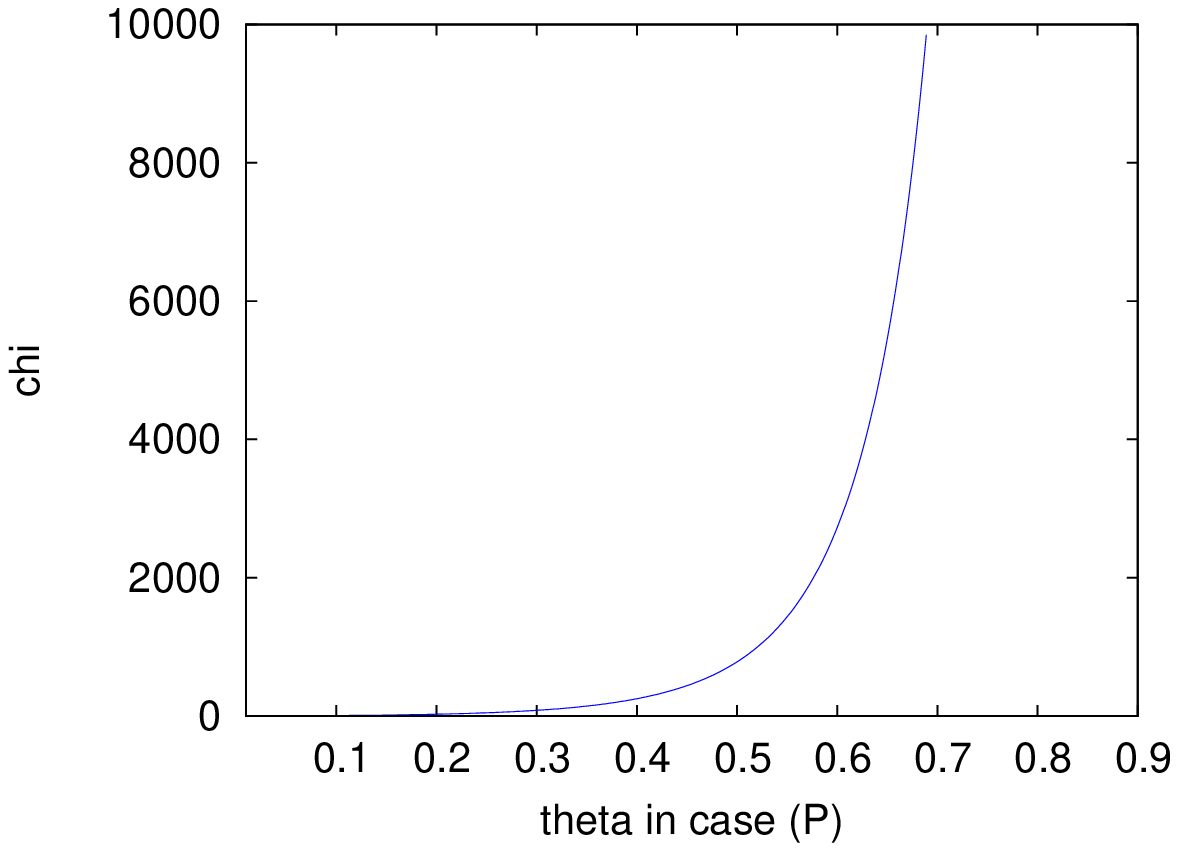}
  \\
  \includegraphics[width=0.49\linewidth]{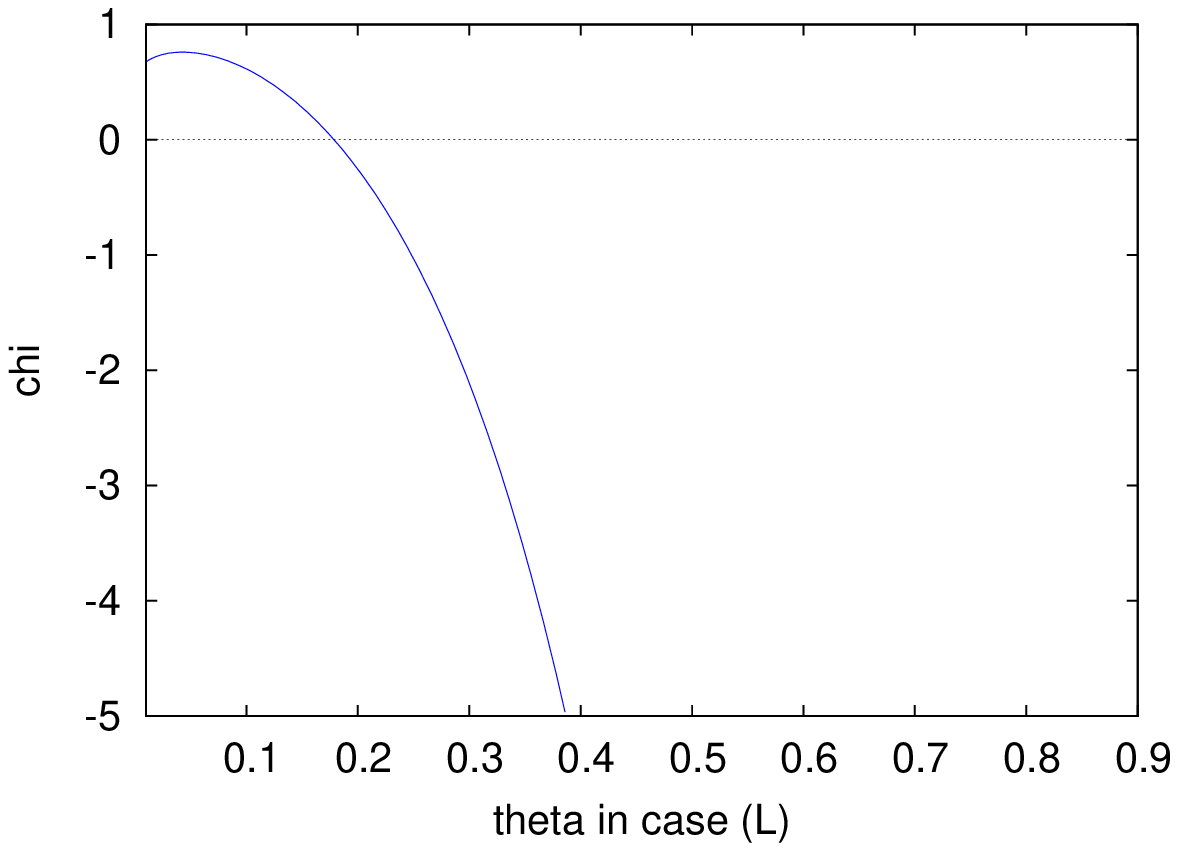}
  \includegraphics[width=0.49\linewidth]{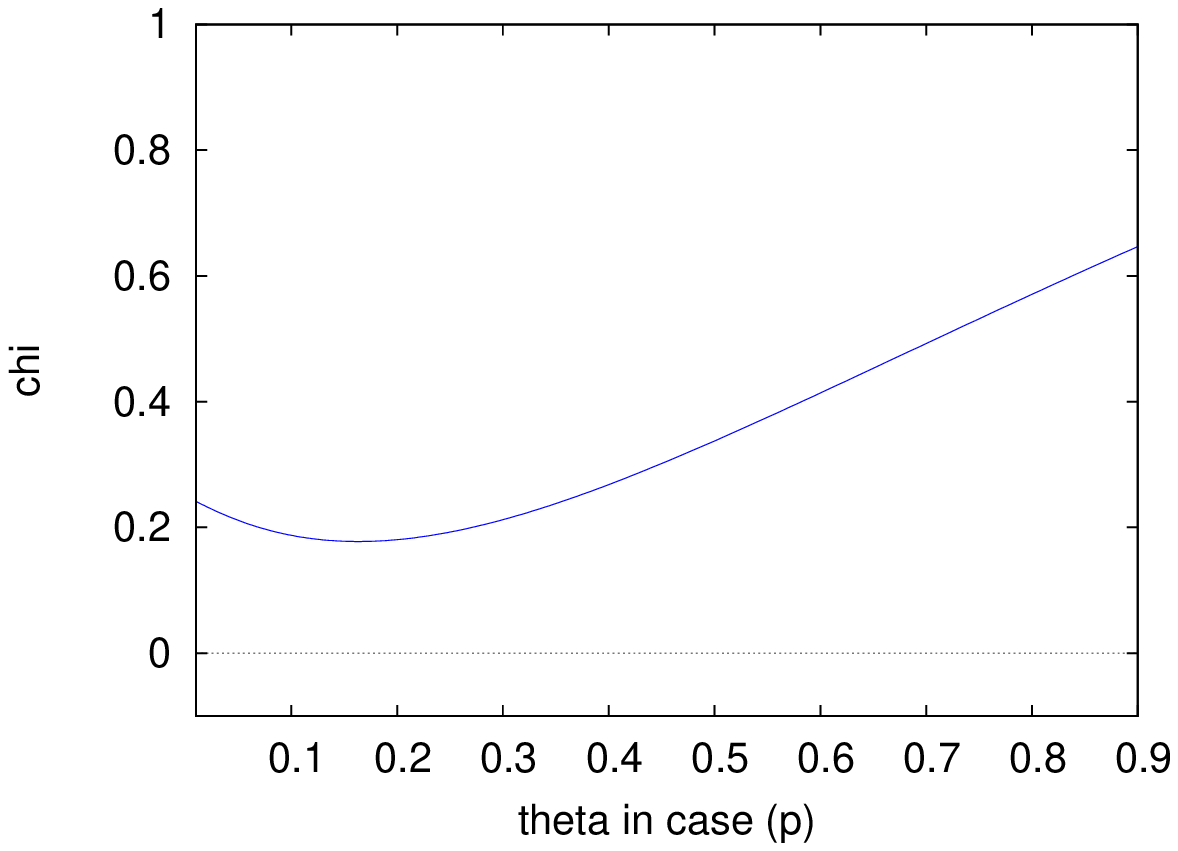}
  \caption{Plots of $\chi$ as a function of $\theta$. Top, left,
    case~\textbf{(S)}; right, case~\textbf{(P)}; bottom, left,
    case~\textbf{(L)}; right, case~\textbf{(p)}. Note that in all four
    cases, $\chi$ attains strictly positive values, showing the
    necessity of the requirement~(\ref{1eq:A0}).}
  \label{1fig:chi}
\end{figure}

The results of these computations are in Figure~\ref{1fig:chi}.  They
show that in all the conditions~(\ref{1eq:Junction}) considered, there
exists a state $u \in A_0$ such that $\chi(\theta) > 0$, showing the
necessity of condition~(\ref{1eq:A0}). However, in case~\textbf{(L)},
it turns out that $\chi$ is negative on an non trivial interval of
values of $\theta$.
If $\bar u$ is chosen in this interval, the wave $\sigma_3$ in the
construction above is not magnified by the consecutive
interactions. The computations leading to the diagrams in
Figure~\ref{1fig:chi} are deferred to Section~\ref{1subs:Conti}.

\Section{Technical Details}
\label{1sec:Tech}

We recall here basic properties of the Euler equations~(\ref{1eq:E1}),
(\ref{1eq:E}). The characteristic speeds and the right eigenvectors
have the expressions
\begin{equation}
  \label{1eigenv}
  \begin{array}{l@{\quad}l@{\quad}l}
    \lambda_{1} 
    = 
    \frac{q}{\rho }-c
    &
    \lambda_{2} 
    = 
    \frac{q}{\rho }
    &
    \lambda_{3} 
    = 
    \frac{q}{\rho }+c
    \\[5pt]
    r_{1}
    =
    \left[
      \!\!\!
      \begin{array}{c}
        -\rho
        \\ 
        \rho c - q
        \\
        qc - E - p
      \end{array}
      \!\!\!
    \right]
    &
    r_{2}
    = 
    \left[ 
      \!\!\!
      \begin{array}{c}
        \rho
        \\ 
        q
        \\
        E+p-\frac{\rho^{2}c^2}{\partial_{e}p}
      \end{array}
      \!\!\!
    \right]
    &
    r_{3} 
    = 
    \left[ 
      \!\!\!
      \begin{array}{c}
        \rho 
        \\ 
        q+\rho c
        \\
        E+p+qc
      \end{array}
      \!\!\!
    \right]
  \end{array}
\end{equation}
whose directions are chosen so that $\nabla \lambda_i \cdot r_i > 0$
for $i=1,2,3$. In the case of an ideal gas, the sound speed $c =
\sqrt{\partial_{\rho}p+\rho^{-2}\,p\,\partial_{e}p}$ becomes
\begin{equation}
  \label{1eq:c}
  c 
  =
  \sqrt{\gamma ( \gamma-1) e} \,.
\end{equation}
The shock and rarefaction curves curves of the first and third family
are:
\begin{eqnarray*}
  S_1(u_o, \sigma) 
  & = &
  \left\{
    \begin{array}{l}
      \rho=-\sigma+\rho_o \\
      v=v_o-\sqrt{-\left(p-p_o\right)\,
        \left(\frac{1}{\rho}-\frac{1}{\rho_o}\right)} \\	  
      e=e_o-\frac{1}{2}\left(p+p_o\right)\,
      \left(\frac{1}{\rho}-\frac{1}{\rho_o}\right)
    \end{array}
  \right.
  \mbox{ for }
  \begin{array}{rcl}
    \sigma & \leq & 0
    \\
    \rho & \geq & \rho_o
    \\
    v & \leq & v_o
    \\
    S & \geq & S_o
  \end{array}
  \\
  S_3(u_o,\sigma)
  & = &
  \left\{
    \begin{array}{l}
      \rho=\sigma+\rho_o \\
      v=v_o-\sqrt{-\left(p-p_o\right)\,
        \left(\frac{1}{\rho}-\frac{1}{\rho_o}\right)} \\	  
      e=e_o-\frac{1}{2}\left(p+p_o\right)\,
      \left(\frac{1}{\rho}-\frac{1}{\rho_o}\right)
    \end{array}
  \right.
  \mbox{ for }
  \begin{array}{rcl}
    \sigma & \leq & 0
    \\
    \rho & \leq & \rho_o
    \\
    v & \leq & v_o
    \\
    S & \leq & S_o
  \end{array}
  \\
  R_1(u_o,\sigma)
  & = & 
  \left\{
    \begin{array}{l}
      \rho=-\sigma+\rho_o \\
      v=v_o-\displaystyle{\int_{p_o}^p [(\rho\,c)(p,S_o)]^{-1}\,dp} \quad\\
      S(\rho,e)=S(\rho_o,e_o)
    \end{array}
  \right.
  \mbox{ for }
  \begin{array}{rcl}
    \sigma& \geq & 0
    \\
    \rho & \leq & \rho_o
    \\
    v & \geq & v_o
    \\
    e & \leq & e_o
  \end{array}
  \\
  R_3(u_o,\sigma)
  & = & 
  \left\{
    \begin{array}{l}
      \rho=\sigma+\rho_o \\
      v=v_o+\displaystyle{\int_{p_o}^p [(\rho\,c)(p,S_o)]^{-1}\,dp} \quad\\
      S(\rho,e)=S(\rho_o,e_o)
    \end{array}
  \right.
  \mbox{ for }
  \begin{array}{rcl}
    \sigma & \geq & 0
    \\
    \rho & \geq & \rho_o
    \\
    v & \geq & v_o
    \\
    e & \geq & e_o
  \end{array}
\end{eqnarray*}
\noindent The 1,2,3-Lax curves have the expressions
\begin{eqnarray*}
L_1(\sigma;\rho_o,q_o,E_o) & = & \left\{
  \begin{array}{lr}
    S_1(\sigma;\rho_o,q_o,E_o), &  \sigma<0 \\
    R_1(\sigma;\rho_o,q_o,E_o), &  \sigma\geq 0 
  \end{array}
\right.
\\
L_2(\sigma;\rho_o,q_o,E_o) & = & \left\{
  \begin{array}{l}
    \rho = \sigma+\rho_o \\
    v = v_o  \\
    p(\rho,e) = p(\rho_o,e_o)
  \end{array}
\right.
\\
L_3(\sigma;\rho_o,q_o,E_o) & = & \left\{
  \begin{array}{lr}
    S_3(\rho;\rho_l,q_l,E_l), &  \sigma<0 \\
    R_3(\sigma;\rho_o,q_o,E_o), &  \sigma\geq 0 
  \end{array}
\right.
\end{eqnarray*}
Their reversed counterparts are
\begin{eqnarray*}
  S_1^{-}(u_o, \sigma) 
  & = &
  \left\{
    \begin{array}{l}
      \rho=\sigma+\rho_o \\
      v=v_o+\sqrt{-\left(p-p_o\right)\,
        \left(\frac{1}{\rho}-\frac{1}{\rho_o}\right)} \\	  
      e=e_o-\frac{1}{2}\left(p+p_o\right)\,
      \left(\frac{1}{\rho}-\frac{1}{\rho_o}\right)
    \end{array}
  \right.
  \mbox{ for }
  \begin{array}{rcl}
    \sigma & \leq & 0
    \\
    \rho & \leq & \rho_o
    \\
    v & \geq & v_o
    \\
    S & \leq & S_o
  \end{array}
  \\
  S_3^{-}(u_o,\sigma)
  & = &
  \left\{
    \begin{array}{l}
      \rho=-\sigma+\rho_o \\
      v=v_o+\sqrt{-\left(p-p_o\right)\,
        \left(\frac{1}{\rho}-\frac{1}{\rho_o}\right)} \\	  
      e=e_o-\frac{1}{2}\left(p+p_o\right)\,
      \left(\frac{1}{\rho}-\frac{1}{\rho_o}\right)
    \end{array}
  \right.
  \mbox{ for }
  \begin{array}{rcl}
    \sigma & \leq & 0
    \\
    \rho & \geq & \rho_o
    \\
    v & \geq & v_o
    \\
    S & \geq & S_o
  \end{array}
  \\
  R_1^{-}(u_o,\sigma)
  & = & 
  \left\{
    \begin{array}{l}
      \rho=\sigma+\rho_o \\
      v=v_o-\displaystyle{\int_{p_o}^p [(\rho\,c)(p,S_o)]^{-1}\,dp} \quad\\
      S(\rho,e)=S(\rho_o,e_o)
    \end{array}
  \right.
  \mbox{ for }
  \begin{array}{rcl}
    \sigma& \geq & 0
    \\
    \rho & \geq & \rho_o
    \\
    v & \leq & v_o
    \\
    e & \geq & e_o
  \end{array}
  \\
  R_3^{-}(u_o,\sigma)
  & = & 
  \left\{
    \begin{array}{l}
      \rho=-\sigma+\rho_o \\
      v=v_o+\displaystyle{\int_{p_o}^p [(\rho\,c)(p,S_o)]^{-1}\,dp}
      \quad \\
      S(\rho,e)=S(\rho_o,e_o)
    \end{array}
  \right.
  \mbox{ for }
  \begin{array}{rcl}
    \sigma & \geq & 0
    \\
    \rho & \leq & \rho_o
    \\
    v & \leq & v_o
    \\
    e & \leq & e_o
  \end{array}
\end{eqnarray*}
and
\begin{eqnarray*}
  L_1^-(\sigma;\rho_o,q_o,E_o) & = & \left\{
    \begin{array}{lr}
      S_1^-(\sigma;\rho_o,q_o,E_o), &  \sigma<0 \\
      R_1^-(\sigma;\rho_o,q_o,E_o), &  \sigma\geq 0 
    \end{array}
  \right.
  \\
  L_2^-(\sigma;\rho_o,q_o,E_o) & = & \left\{
    \begin{array}{l}
      \rho=-\sigma+\rho_o \\
      v=v_o  \\
      p(\rho,e)=p(\rho_o,e_o)
    \end{array}
  \right.
  \\
  L_3^-(\sigma;\rho_o,q_o,E_o) & = & \left\{
    \begin{array}{ll}
      S_3^-(\sigma;\rho_o,q_o,E_o), &  \sigma<0 \\
      R_3^-(\sigma;\rho_o,q_o,E_o), &  \sigma\geq 0 \,.
    \end{array}
  \right.
\end{eqnarray*}

In the $(\rho, q, e)$ space, for a perfect ideal gas, the tangent
vectors to the Lax curves are:
\begin{eqnarray}
  \label{1tangent vectors}
  \tilde r_1
  =
  \left[\!\!
    \begin{array}{c}
      -1
      \\
      -\frac{q}{\rho} -\sqrt{\gamma ( \gamma-1) e}
      \\
      -(\gamma-1) \frac{e}{\rho}
    \end{array}
  \!\!\right]
  \,,\  
  \tilde r_2
  =
  \left[\!\!
    \begin{array}{c}
      1
      \\
      \frac{q}{\rho}
      \\
      -\frac{e}{\rho}
    \end{array}
  \!\!\right]
  \,,\ 
  \tilde r_3
  =
  \left[\!\!
    \begin{array}{c}
      1
      \\
      \frac{q}{\rho} -\sqrt{\gamma ( \gamma-1) e}
      \\
      (\gamma-1) \frac{e}{\rho}
    \end{array}
  \!\!\right].\;
\end{eqnarray}

\subsection{Proofs of Section~\ref{1sec:General}}

The following result will be of use in the proof of
Proposition~\ref{1thm:CauchyProblem}.
\begin{proposition}
  \label{1prop:Lax curves}
  Let $\sigma_{i}\mapsto L_{i}(u_{0},\sigma_{i})$ be the $i$-th Lax
  curve and $\sigma_{i}\mapsto L_{i}^{-}(u_{0},\sigma_{i})$ be the
  reversed $i$-th Lax curve through $u_{0}$, for $i=1,2,3$. The
  following equalities hold:
  \begin{displaymath}
    \frac{\partial L_1}{\partial\sigma_{1}}_{|\sigma_{1}=0}= \left(
      \begin{array}{c}
        1 \\
        \lambda_1(u_o) \\
        \displaystyle{\frac{E_o+p_o}{\rho_o}
          -\frac{q_o}{\rho_o}\,c_o}
      \end{array}
    \right),\ \frac{\partial L_2}{\partial\sigma_{2}}_{|\sigma_{2}=0}=
    \left(
      \begin{array}{c}
        1 \\
        \lambda_2(u_o) 
        \\
        \displaystyle
        {\frac{E_o+p_o}{\rho_o}-\frac{\rho_o\,c_o^2}{\partial_e p_o}}
      \end{array}
    \right),
  \end{displaymath}
  \begin{displaymath}
    \frac{\partial L_3}{\partial\sigma_{3}}_{|\sigma_{3}=0}= \left(
      \begin{array}{c}
        1 \\
        \lambda_3(u_o) \\
        \displaystyle{\frac{E_o+p_o}{\rho_o}+\frac{q_o}{\rho_o}\,c_o}
      \end{array}
    \right),
  \end{displaymath}
  \begin{displaymath}
    \mbox{ for }i=1,2,3\qquad \frac{\partial
      L_i^-}{\partial\sigma_{i}}_{|\sigma_{i}=0}= -\frac{\partial
      L_i}{\partial\sigma_{i}}_{|\sigma_{i}=0},
  \end{displaymath}
  \begin{displaymath}
    \frac{\partial L_i}{\partial \rho_o}_{|\sigma_{i}=0}= \left(
      \begin{array}{c}
        1 \\
        0 \\
        0
      \end{array}
    \right),\quad \frac{\partial L_i}{\partial q_o}_{|\sigma_{i}=0}=
    \left(
      \begin{array}{c}
        0 \\
        1 \\
        0
      \end{array}
    \right),\quad \frac{\partial L_i}{\partial E_o}_{|\sigma_{i}=0}=
    \left(
      \begin{array}{c}
        0 \\
        0 \\
        1
      \end{array}
    \right),
  \end{displaymath}
  \begin{displaymath}
    +    \frac{\partial L_i^-}{\partial \rho_o}_{|\sigma_{i}=0}= \frac{\partial
      L_i}{\partial \rho_o}_{|\sigma_{i}=0},\quad \frac{\partial
      L_i^-}{\partial q_o}_{|\sigma_{i}=0}= \frac{\partial L_i}{\partial
      q_o}_{|\sigma_{i}=0},\quad \frac{\partial L_i^-}{\partial
      E_o}_{|\sigma_{i}=0}= \frac{\partial L_i}{\partial
      E_o}_{|\sigma_{i}=0}.
  \end{displaymath}
\end{proposition}
The proof is immediate and, hence, omitted.

\begin{proofof}{Theorem~\ref{1thm:CauchyProblem}}
  Following~\cite[Proposition~4.2]{ColomboGuerraHertySachers}, the
  $3\times 3$ system~(\ref{1eq:E}) defined for $x \in \reali$ can be
  rewritten as the following $6\times 6$ system defined for $x \in
  \reali^+$:
  \begin{equation}
    \label{1eq:UF}
    \left\{
      \begin{array}{l@{\qquad} rcl}
        \partial_t U + \partial_x \mathcal{F}(U) = 0
        &
        (t,x) & \in & \reali^+ \times \reali^+
        \\
        b \left( U(t,0+) \right) = 0
        &
        t & \in & \reali^+
      \end{array}
    \right.
  \end{equation}
  the relations between $U$ and $u = (\rho, q, E)$, between
  $\mathcal{F}$ and the flow in~(\ref{1eq:E}) being
  \begin{displaymath}
    U(t,x)
    =
    \left[
      \begin{array}{c}
        \rho(t, -x) \\ q(t, -x) \\ E(t, -x) \\
        \rho(t, x) \\ q(t, x) \\ E(t, x)
      \end{array}
    \right]
    \quad\mbox{ and } \quad
    \mathcal{F}(U)
    =
    \left[
      \begin{array}{c}
        U_2 \\ P(U_1,U_2,U_3) \\ F(U_1,U_2,U_3) \\ 
        U_5 \\ P(U_4,U_5,U_6) \\ F(U_4,U_5,U_6)
      \end{array}
    \right]
  \end{displaymath}
  with $x \in \reali^+$ and $E,P,F$ defined in~(\ref{1eq:EPF}); whereas
  the boundary condition in~(\ref{1eq:UF}) is related
  to~(\ref{1eq:Junction}) by
  \begin{displaymath}
    b(U) = \Psi \left(a^-, (U_1, U_2, U;_3); a^+, (U_4, U_5, U_6) \right)
  \end{displaymath}
  for fixed sections $a^-$ and $a^+$.

  The thesis now follows
  from~\cite[Theorem~2.2]{ColomboGuerra6}. Indeed, the
  assumptions~\textbf{($\mathbf{\gamma}$)}, \textbf{(b)}
  and~\textbf{(f)} therein are here satisfied. More precisely,
  condition~\textbf{($\mathbf{\gamma}$)} follows from the
  choice~(\ref{1subsonic}) of the subsonic region $A_0$. Simple
  computations show that condition~\textbf{(b)} reduces to
  \begin{displaymath}
    \begin{array}{l}
      \det \left[ 
        D_{u^{-}}\Psi \cdot 
        \frac{\partial L_1}{\partial\sigma_{1}}_{|\sigma_{1}=0}
        \quad
        D_{u^{+}}\Psi \cdot 
        \frac{\partial L_2^{-}}{\partial\sigma_{2}}_{|\sigma_{2}=0}
        \quad
        D_{u^{+}}\Psi \cdot 
        \frac{\partial L_2^{-}}{\partial u^{+}}_{|\sigma_{2}=0} \cdot 
        \frac{\partial L_3^{-}}{\partial\sigma_{3}}_{|\sigma_{3}=0}
      \right]    
      \\
      =
      \det \left[ 
        D_{u^{-}}\Psi\cdot r_{1}(\bar u)
        \quad
        -D_{u^{+}}\Psi \cdot r_{2}(\bar u)
        \quad
        -D_{u^{+}}\Psi \cdot r_{3}(\bar u)
      \right]
      \\
      =
      \det \left[ 
        D_{u^{-}}\Psi\cdot r_{1}(\bar u)
        \quad
        D_{u^{+}}\Psi \cdot r_{2}(\bar u)
        \quad
        D_{u^{+}}\Psi \cdot r_{3}(\bar u)
      \right],
    \end{array}
  \end{displaymath}
  which is non zero for assumption if $\bar u\in A_{0}$ and $\bar
  a > 0$. Condition~\textbf{(f)} needs more care. Indeed,
  system~(\ref{1eq:UF}) is \emph{not} hyperbolic, for it is obtained
  gluing two copies of the Euler equations~(\ref{1eq:E}). Nevertheless,
  the two systems are coupled only through the boundary condition,
  hence the whole wave front tracking procedure in the proof
  of~\cite[Theorem~2.2]{ColomboGuerra6} applies, see
  also~\cite[Proposition~4.5]{ColomboGuerraHertySachers}.
\end{proofof}

\begin{proofof}{Proposition~\ref{1prop:otherpsi}}
  It is immediate to check that each of the coupling
  conditions~(\ref{1eq:psismooth}), (\ref{1eq:psiCri}), (\ref{1eq:psiL}),
  (\ref{1eq:psip}) satisfies the
  requirements~(\textbf{$\mathbf{\Psi}$0})
  and~(\textbf{$\mathbf{\Psi}$1}).

  To prove that~(\textbf{$\mathbf{\Psi}$2}) is satisfied, we use an
  \emph{ad hoc} argument for condition~\textbf{(S)}. In all the other
  cases, note that the function $\Psi$ admits the representation
  $\Psi(a^-, u^-; a^+, u^+) = \psi(a^-, u^-) -
  \psi(a^+,u^+)$. Therefore, (\textbf{$\mathbf{\Psi}$2}) trivially
  holds.

  We prove below~(\ref{1eq:det}) in each case separately. Note however
  that for any of the considered choices of $\Psi$,
  \begin{eqnarray}
    \label {condition}
    D_{u^{+}}\Psi(a^{-},u^{-},a^{+},u^{+})_{|u=\bar u,\,a=\bar a}=-D_{u^{-}}\Psi(a^{-},u^{-},a^{+},u^{+})_{|u=\bar u,\,a=\bar a}
  \end{eqnarray}
  so that~(\ref{1eq:det}) reduces to
  \begin{eqnarray*}
    & &
    \det \left[ 
      D_{u^{-}}\Psi\cdot r_{1}(\bar u)
      \quad
      D_{u^{+}}\Psi \cdot r_{2}(\bar u)
      \quad
      D_{u^{+}}\Psi \cdot r_{3}(\bar u)
    \right]    
    \\
    & = &
    -\det D_{u^{+}}\Psi \cdot
    \det \left[ 
      r_{1}(\bar u)
      \quad
      r_{2}(\bar u)
      \quad
      r_{3}(\bar u)
    \right] \,.
  \end{eqnarray*}
  Thus, it is sufficient to prove that $\det
  D_{u^{+}}\Psi(a^{-},u^{-},a^{+},u^{+})_{|u=\bar u,\,a=\bar a}\neq
  0$.

  \paragraph{\textbf{(S)-solutions}}

  To prove that the coupling condition~(\ref{1eq:psismooth})
  satisfies~(\textbf{$\mathbf{\Psi}$2}), simply use the additivity of
  the integral and the uniqueness of the solution to the Cauchy
  problem for the ordinary differential
  equation~(\ref{1eq:Stationary}).

  Next, we have
  \begin{displaymath}
    D_{u}
    \left( 
      \int_{-X}^X p \left( \mathcal{R}^{a}(x), \mathcal{E}^a(x) \right) 
      \, a'(x)
      \, \mathrm{d}x
    \right)_{|u=\bar u,\,a=\bar a} = 0 \,,
  \end{displaymath}
  since $a'(x)=0$ for all $x$, because $a^- = a^+ = \bar a$.  Thus,
  $\Psi$ in~(\ref{1eq:psismooth}) satisfies
  \begin{displaymath}
    \begin{array}{l}
      D_{u^{+}}\Psi(a^{-},u^{-},a^{+},u^{+})_{|u=\bar u,\,a=\bar a}
      \\
      =  
      \bar a^{3}\det\left[  
        \!\!
        \begin{array}{ccc}
          0  & 1 & 0
          \\ 
          -\frac{\bar q^{2}}{\bar \rho^{2}}+{\partial_\rho \bar
            p}+\frac{{\partial_e \bar p}}{\bar \rho}\left( \frac{\bar
              q^{2}}{\bar \rho^{2}}-\frac{\bar E}{\bar\rho}\right) 
          &
          \frac{\bar q}{\bar\rho}\left( 2-\frac{{\partial_e \bar p}}{\bar
              \rho}\right)
          &  
          \frac{{\partial_e \bar p}}{\bar \rho}
          \\ 
          -\frac{\bar q}{\bar\rho}\left( {\partial_\rho \bar
              p}+\frac{{\partial_e \bar p}}{\bar \rho}\left( \frac{\bar
                q^{2}}{\bar \rho^{2}}-\frac{\bar E}{\bar\rho}\right)
            -\frac{\bar E+\bar p}{\bar\rho}\right)
          &
          \frac{\bar E+\bar p}{\bar\rho}-\frac{{\partial_e \bar p}}{\bar
            \rho}\frac{\bar q^{2}}{\bar \rho^{2}}
          &
          -\frac{\bar q}{\bar\rho}\left( \frac{{\partial_e \bar p}}{\bar \rho}+1\right) 
        \end{array}
        \!\!
      \right]
\!\!
      \\
      =  
      -\bar a^{3} \, 
      \lambda_{1}(\bar u) \, \lambda_{2}(\bar u) \, \lambda_{3}(\bar u),
    \end{array}
  \end{displaymath}
  which is non zero if $\bar u\in A_{0}$.

  \paragraph{\textbf{(P)-solutions}}

  Concerning condition~(\ref{1eq:psiCri}), we have
  \begin{displaymath}
    \begin{array}{l}
      D_{u^{+}}\Psi(a^{-},u^{-},a^{+},u^{+})_{|u=\bar u,\,a=\bar a}
      \\
      =
      \det\left[  
\!\!
        \begin{array}{ccc}
          0 & \bar a & 0
          \\ 
          -\frac{\bar q^{2}}{\bar \rho^{2}}+{\partial_\rho \bar
            p}+\frac{{\partial_e \bar p}}{\bar \rho}\left( \frac{\bar
              q^{2}}{\bar \rho^{2}}-\frac{\bar E}{\bar\rho}\right) 
          &\frac{\bar q}{\bar\rho}\left( 2-\frac{{\partial_e \bar p}}{\bar
              \rho}\right)  
          & \frac{{\partial_e \bar p}}{\bar \rho}
          \\ 
          -\bar a\frac{\bar q}{\bar\rho}\left( {\partial_\rho \bar
              p}+\frac{{\partial_e \bar p}}{\bar \rho}\left( \frac{\bar
                q^{2}}{\bar \rho^{2}}-\frac{\bar E}{\bar\rho}\right)
            -\frac{\bar E+\bar p}{\bar\rho}\right) 
          & \bar a\frac{\bar E+\bar p}{\bar\rho}-\bar a\frac{{\partial_e \bar
              p}}{\bar \rho}\frac{\bar q^{2}}{\bar \rho^{2}}
          & -\bar a\frac{\bar q}{\bar\rho}\left( \frac{{\partial_e \bar p}}{\bar \rho}+1\right) 
        \end{array}%
      \!\!
\right]
\!\!
      \\
      = 
      -\bar a^{2}\lambda_{1}(\bar u)\lambda_{2}(\bar u)\lambda_{3}(\bar u),
    \end{array}
  \end{displaymath}
  which is non zero if $\bar u\in A_{0}$.

  \paragraph{\textbf{(L)-solution}} For condition~(\ref{1eq:psiL}) the
  computations very similar to the above case:
  \begin{eqnarray*}
    D_{u^{+}}\Psi(a^{-},u^{-},a^{+},u^{+})_{|u=\bar u,\,a=\bar a}=  - \bar a^{3}\lambda_{1}(\bar u)\lambda_{2}(\bar u)\lambda_{3}(\bar u),
  \end{eqnarray*}
  which is non zero if $\bar u\in A_{0}$.
  
  \paragraph{\textbf{(p)-solution}} Finally, concerning
  condition~(\ref{1eq:psip}),
  \begin{displaymath}
    \begin{array}{l}
      D_{u^{+}}\Psi(a^{-},u^{-},a^{+},u^{+})_{|u=\bar u,\,a=\bar a}
      \\
      =
      \det\left[  
        \!\!
        \begin{array}{ccc}
          0 & \bar a & 0\\ 
          {\partial_\rho \bar p}+\frac{{\partial_e \bar p}}{\bar
            \rho}\left( \frac{\bar q^{2}}{\bar \rho^{2}}-\frac{\bar
              E}{\bar\rho}\right)  &-\frac{\bar
            q}{\bar\rho^{2}}{\partial_e \bar p} & \frac{{\partial_e \bar
              p}}{\bar \rho}
          \\ 
          -\bar a\frac{\bar q}{\bar\rho}\left( {\partial_\rho \bar p}+\frac{{\partial_e \bar p}}{\bar \rho}\left( \frac{\bar q^{2}}{\bar \rho^{2}}-\frac{\bar E}{\bar\rho}\right) -\frac{\bar E+\bar p}{\bar\rho}\right)  & \bar a\frac{\bar E+\bar p}{\bar\rho}-\bar a\frac{{\partial_e \bar p}}{\bar \rho}\frac{\bar q^{2}}{\bar \rho^{2}} & -\bar a\frac{\bar q}{\bar\rho}\left( \frac{{\partial_e \bar p}}{\bar \rho}+1\right) 
        \end{array}%
        \!\!
      \right]
      \!\!
      \\
      =
      \bar a^{2}\lambda_{2}(\bar u)\left( c^{2}+\lambda_{2}^{2}(\bar u)\frac{{\partial_e \bar p}}{\bar \rho}\right) ,
    \end{array}
  \end{displaymath}
  which is non zero if $\bar u\in A_{0}$ and if the fluid is perfect,
  i.e.~(\ref{1eq:perfect}) holds.
\end{proofof}

\subsection{Computation of $\chi$ in~(\ref{1eq:chi})}
\label{1subs:Conti}

\paragraph{The Case of Condition~(S)}

Let $\Psi$ be defined in~(\ref{1eq:psismooth}) and set
\begin{displaymath}
  \Sigma(a^-, a^+, u)
  =
  \int_{-X}^X p\left( \mathcal{R}^a(x), \mathcal{E}^a(x) \right) \,
  a'(x) d x
\end{displaymath}
where the functions $\mathcal{R}^a, \mathcal{E}^a$ have the same
meaning as in~(\ref{1eq:psismooth}). A perturbative method allows to
compute the solution to~(\ref{1eq:Stationary}) with a second order
accuracy in $(\Delta a)/a$. Then, long elementary computations allow
to get explicitly the terms $H$ and $G$ in~(\ref{1eq:TT}) of the second
order expansion of $T$:
\begin{eqnarray*}
  H (\rho, q, e)
  & = &
  \left[
    \begin{array}{l}
      \displaystyle
      -{{\,\theta^3-4\,\theta^2+5\,\theta-2}\over{\,
          \theta^3-3\,\theta^2+3\,\theta-1}}
      \; \rho
      \\
      \displaystyle
      -q
      \\
      \displaystyle
      -{{2\left( -\theta^3+2\,\theta^2-\theta\right) }\over{3\left(\theta^3
            -3\,\theta^2+3\,\theta-1\right) }}
      \; e
    \end{array}
  \right]
  \\
  G(\rho, q, e)
  & = &
  \left[
    \begin{array}{l}
      \displaystyle
      -{{4\left(\theta^3-2\,\theta^2\right) }\over{3\left( \theta^3-3\,\theta^2
            +3\,\theta-1\right)}}
      \; \rho
      \\
      q
      \\
      \displaystyle
      -{{70\,\theta^4-257\,\theta^3+342\,\theta^2-207\,\theta
          +36}\over{18\left( \theta^3-3\,\theta^2+3\,\theta-1\right) }}
      \; e
    \end{array}
  \right]\,.
\end{eqnarray*}
Moreover, the coefficients $f_1$, $f_2$ in~(\ref{1first}) read
\begin{eqnarray*}
  f_1(\theta)
  & = &
  -{{-3\,\theta+\left(\theta-3\right)\,\sqrt{\theta}-3
    }\over{6\,\sqrt{\theta}\,\left(\theta-1\right)-6\left(\theta-1\right) 
    }}
  \\
  f_2(\theta)
  & = &
  {{\sqrt{\theta}\,\left(126\,\theta^4-505\,\theta^3+758\,
        \theta^2-489\,\theta+270\right)\over{72\left( \sqrt{\theta}\,\left(       \theta^3-3\,\theta^2+3\,\theta-1\right)-
          \theta^3+3\,\theta^2-3\,\theta+1\right) }}}
  \\
  & &
  +{{42\,\theta^4-183\,\theta
      ^3+278\,\theta^2+33\,\theta+54}\over{72\left( \sqrt{\theta}\,\left(       \theta^3-3\,\theta^2+3\,\theta-1\right)-
        \theta^3+3\,\theta^2-3\,\theta+1\right)}}\,.
\end{eqnarray*}
Next, $\chi$ is given by {\footnotesize
  \begin{displaymath}
    \chi
    =
    {{\sqrt{\theta}\,\left(126\,\theta^4-506\,\theta^3+773\,
          \theta^2-480\,\theta+279\right)+42\,\theta^4-174\,\theta
        ^3+311\,\theta^2+96\,\theta+45}\over{36\left( \sqrt{\theta}\,\left(
            \theta^3-3\,\theta^2+3\,\theta-1\right)-
          \theta^3+3\,\theta^2-3\,\theta+1\right) }}\,.
  \end{displaymath}
}

\paragraph{The Case of Condition~\textbf{(P)}}

Let $\Psi$ be defined in~(\ref{1eq:psiCri}). With reference
to~(\ref{1eq:TT}), we show below explicitly the terms $H$ and $G$
in~(\ref{1eq:TT}) of the second order expansion of $T$,
\begin{eqnarray*}
  H (\rho, q, e)
  & = &
  \left[
    \begin{array}{l}
      \displaystyle
      {{8\left( -\theta^3+2\,\theta^2-\theta\right) }\over{3\left( 
            \theta^3-3\,\theta^2+3\,\theta-1\right) }} 
      \; \rho
      \\
      \displaystyle
      -q
      \\
      \displaystyle
      -{{2\left(5\, \theta^4-7\,\theta^3-\theta^2+3\,\theta\right) 
        }\over{9\left( \theta^3-3\,\theta^2+3\,\theta-1\right) }}
      \; e
    \end{array}
  \right]
  \\
  G(\rho, q, e)
  & = &
  \left[
    \begin{array}{l}
      \displaystyle
      {{64\left( \theta^3+3\,\theta^2\right) }\over{27\left( \theta^3-3\,
            \theta^2+3\,\theta-1\right) }}
      \; \rho
      \\
      q
      \\
      \displaystyle
      -{{565\,\theta^4-1599\,\theta^3+927\,\theta^2-405\,
          \theta}\over{81\left( \theta^3-3\,\theta^2+3\,\theta-1\right) }}
      \; e
    \end{array}
  \right]\,.
\end{eqnarray*}
Moreover, the coefficients $f_1$, $f_2$ in~(\ref{1first}) read
\begin{eqnarray*}
  f_1(\theta)
  & = &
  {{\sqrt{\theta}\,\left(9\,\theta^2+2\,\theta-27\right)+3\,
      \theta^2-42\,\theta-9}\over{18\sqrt{\theta}\,\left(
        \theta-1\right)-18\left( \theta-1\right) }}
  \\
  f_2(\theta)
  & = &
  {{\sqrt{\theta}\,\left(154\,\theta^5+931\,\theta^4-4416\,
        \theta^3+6570\,\theta^2+990\,\theta+891\right)}\over{324\left( \sqrt{\theta}\,\left(\theta^3-3\,
          \theta^2+3\,\theta-1\right)-\theta^3+3\,
        \theta^2-3\,\theta+1\right) }}
  \\
  & &
  +{{86\,
      \theta^5-311\,\theta^4-752\,\theta^3+7038\,\theta^2+1026
      \,\theta+81}\over{324\left( \sqrt{\theta}\,\left(\theta^3-3\,
          \theta^2+3\,\theta-1\right)-\theta^3+3\,
        \theta^2-3\,\theta+1\right) }}.
\end{eqnarray*}
Next, $\chi$ is given by
\begin{eqnarray*}
  \chi
  & = &
  {{\sqrt{\theta}\,\left(407\,\theta^5+1931\,\theta^4-7858\,
        \theta^3+14766\,\theta^2+1179\,\theta+1863\right)
    }\over{324\left( \sqrt{\theta}\,\left(\theta^3-3\,
          \theta^2+3\,\theta-1\right)-\theta^3+3\,
        \theta^2-3\,\theta+1\right)}}
  \\
  & &
  + 
  {{-23\,
      \theta^5+141\,\theta^4+2002\,\theta^3+15714\,\theta^2+2565
      \,\theta+81}\over{324\left( \sqrt{\theta}\,\left(\theta^3-3\,
          \theta^2+3\,\theta-1\right)-\theta^3+3\,
        \theta^2-3\,\theta+1\right)}}\,. 
\end{eqnarray*}

\paragraph{The Case of Condition~\textbf{(L)}}

Let $\Psi$ be defined in~(\ref{1eq:psiL}). Then,
\begin{displaymath}
  H(\rho, q, e)
  =
  \left[
    \begin{array}{l}
      -\rho
      \\
      -q
      \\
      0
    \end{array}
  \right]
  \qquad \mbox{ and }\qquad
  G(\rho, q, e)
  = 
  \left[
    \begin{array}{l}
      \displaystyle
      -{{4\,\theta}\over{3\left( \theta-1\right) }} \; \rho
      \\
      q
      \\
      \displaystyle
      -{{35\,\theta^2-9\left( 4\,\theta-1\right) }\over{9\left( \theta-1\right) }} \; e
    \end{array}
  \right].
\end{displaymath}
The coefficients $f_1$, $f_2$ in~(\ref{1first}) read
\begin{eqnarray*}
  f_1(\theta)
  & = &
  0
  \\
  f_2(\theta)
  & = &
  {{\sqrt{\theta}\,\left(63\,\theta^2-106\,\theta+27\right)+21
      \,\theta^2-78\,\theta+9}\over{36\left( \sqrt{\theta}\,\left(
          \theta-1\right)-\theta+1\right) }},
\end{eqnarray*}
so that $\chi$ is
\begin{displaymath}
  \chi
  =
  {{\sqrt{\theta}\,\left(63\,\theta^2-106\,\theta+27\right)+21
      \,\theta^2-78\,\theta+9}\over{18\left( \sqrt{\theta}\,\left(
          \theta-1\right)-\theta+1\right) }}\,.
\end{displaymath}

\paragraph{The Case of Condition~\textbf{(p)}}

Let $\Psi$ be defined in~(\ref{1eq:psip}).  With reference
to~(\ref{1eq:TT}),
\begin{eqnarray*}
  H(\rho, q, e)
  & = &
  \left[
    \begin{array}{l}
      \displaystyle
      -{{2\left( 4\,\theta^3+12\,\theta^2+9\,\theta\right) }\over{4\left( 2\,\theta^3
            +9\,\theta^2\right) +27\left( 2\,\theta+1\right) }}
      \; \rho
      \\
      -q
      \\
      \displaystyle
      {{2\left( 4\,\theta^3+12\,\theta^2+9\,\theta\right) }\over{4\left( 2\,\theta^3
            +9\,\theta^2\right) +27\left( 2\,\theta+1\right)}} \; e
    \end{array}
  \right]
  \\
  G(\rho, q, e)
  & = &
  \left[
    \begin{array}{l}
      \displaystyle
      -{{4\left( \theta^3+3\,\theta^2\right) }\over{4\left( 2\,\theta^3
            +9\,\theta^2\right) +27\left( 2\,\theta+1\right)}}
      \; \rho
      \\
      q
      \\
      \displaystyle
      {{12\left( \theta^3+2\,\theta^2\right) }\over{4\left( 2\,\theta^3
            +9\,\theta^2\right) +27\left( 2\,\theta+1\right)}}
      \; e
    \end{array}
  \right],
\end{eqnarray*}
with $f_1$ and $f_2$ given by
\begin{eqnarray*}
  f_1(\theta)
  & = &
  {{-2\,\theta^2+4\,\theta^{{{3}\over{2}}}+3\,\theta-9
    }\over{2\left( 4\,\theta^2+12\,\theta+9\right) }}
  \\
  f_2(\theta)
  & = &
  {{32\,\theta^4+8\,\sqrt{\theta}\,\left(4\,\theta^3+9\,
        \theta^2-9\,\theta\right)+316\,\theta^3+558\,\theta^2+216
      \,\theta+81}\over{6\left( 16\,\theta^4+96\,\theta^3+216\,\theta
        ^2+216\,\theta+81\right) }}\,,
\end{eqnarray*}
so that
\begin{displaymath}
  \chi
  =
  {{60\,\theta^4+96\,\sqrt{\theta}\,\left(\theta^3+
        \theta^2-3\,\theta\right)+700\,\theta^3+1107\,\theta^2
      -54\,\theta+81}\over{6\left( 16\,\theta^4+96\,\theta^3+216\,\theta
        ^2+216\,\theta+81\right)}}.
\end{displaymath}

{\footnotesize
  \begin{landscape}
    \ \vspace{-5\baselineskip}
    \begin{eqnarray*}
      \theta^+
      & = &
      \theta
      -
      {{\sqrt{\theta}\,\left(\left(\sigma_3^-+6\right)\,\theta^2
            +18\,\theta-9\,\sigma_3^-\right)+2\left(3-2\,\sigma_3^-
          \right)\,\theta^2+6\left(2\,\sigma_3^-+3\right)\,\theta
        }\over{9\sqrt{\theta}\,\left(\theta-1\right)+9\left( \theta-1\right) 
        }}
      \frac{\Delta a}{a}
      \\
      & - &
      {{\sqrt{\theta}\,\left(14\left(11\,\sigma_3^--30\right)\,
            \theta^5+\left(990-301\,\sigma_3^-\right)\,\theta^4+3\left(25
              \,\sigma_3^--236\right)\,\theta^3+\left(111\,\sigma_3^--18
            \right)\,\theta^2+45\left(12-\,\sigma_3^-\right)\,\theta-378
            \,\sigma_3^-\right)}\over{108\left( \sqrt{\theta}\,\left(\theta^3-3\,\theta^2+3
              \,\theta-1\right)+\theta^3-3\,\theta^2+3\,
            \theta-1\right) }}
      \left( \frac{\Delta a}{a} \right)^2
      \\
      & &
      -
      {{4\left(217\,\sigma_3^--105\right)\,
          \theta^5+2\left(495-1489\,\sigma_3^-\right)\,\theta^4+4\left(928
            \,\sigma_3^--177\right)\,\theta^3-2\left(1077\,\sigma_3^--9
          \right)\,\theta^2+24\left(26\,\sigma_3^-+15\right)\,\theta
        }\over{108\left( \sqrt{\theta}\,\left(\theta^3-3\,\theta^2+3
              \,\theta-1\right)+\theta^3-3\,\theta^2+3\,
            \theta-1\right)}}
      \left( \frac{\Delta a}{a} \right)^2
      \\
      \theta^+
      & = &
      \theta
      -    {{\sqrt{\theta}\,\left(\left(11\,\sigma_3^--30\right)\,
            \theta^3+\left(-19\,\sigma_3^--108\right)\,\theta^2+9\left(5
              \,\sigma_3^--6\right)\,\theta+27\,\sigma_3^-\right)+
          2\left(31
            \,\sigma_3^--15\right)\,\theta^3-36\left(\sigma_3^-+3
          \right)\,\theta^2-18\left(5\,\sigma_3^-+3\right)\,\theta
        }\over{27\sqrt{\theta}\,\left(\theta-1\right)+27\left( \theta
            -1\right) }}
      \frac{\Delta a}{a}
      \\
      & - &
      {{\sqrt{\theta}\,\left(2\left(233\,\sigma_3^--300\right)\,
            \theta^6+\left(1279\,\sigma_3^--5310\right)\,\theta^5+
            \left(5400-2543\,\sigma_3^-\right)\,\theta^4+6\left(677\,
              \sigma_3^-+1242\right)\,\theta^3+36\left(109-297\,
              \sigma_3^-\right)\,\theta^2+729\left(2-5\,\sigma_3^-
            \right)\,\theta-1215\,\sigma_3^-\right)}\over{486\left( \sqrt{\theta}\,\left(
              \theta^3-3\,\theta^2+3\,\theta-1\right)+
            \theta^3-3\,\theta^2+3\,\theta-1\right) }}
      \left(  \frac{\Delta a}{a} \right)^2
      \\
      & &
      -
      {{4\left(601\,
            \sigma_3^--150\right)\,\theta^6+2\left(3521\,\sigma_3^--2655
          \right)\,\theta^5+8\left(225-3814\,\sigma_3^-\right)\,
          \theta^4+36\left(655\,\sigma_3^-+207\right)\,\theta^3+
          108\left(53\,\sigma_3^-+36\right)\,\theta^2+162\left(9\,
            \sigma_3^-+25\right)\,\theta}\over{486\left( \sqrt{\theta}\,\left(
              \theta^3-3\,\theta^2+3\,\theta-1\right)+
            \theta^3-3\,\theta^2+3\,\theta-1\right)}}
      \left(  \frac{\Delta a}{a} \right)^2
      \\
      \theta^+
      & = &
      \theta
      -
      {{\sqrt{\theta}\,\left(\left(77\,\sigma_3^--210\right)\,
            \theta^3+\left(-13\,\sigma_3^--36\right)\,\theta^2+27\left(
              \sigma_3^-+2\right)\,\theta-27\,\sigma_3^-\right)
          2+\left(217
            \,\sigma_3^--105\right)\,\theta^3-4\left(147\,\sigma_3^-+9
          \right)\,\theta^2+18 \left(5\,\sigma_3^-+3\right)\,\theta
        }\over{54\sqrt{\theta}\,\left(\theta-1\right)+54\left( \theta
            -1\right) }}
      \left( \frac{\Delta a}{a} \right)^2
      \\
      \theta^+
      & = &
      \theta
      +
      {{-2\left({\sigma_3^-}+6\right)\,\theta^3+\sqrt{\theta}\,
          \left(10\,{\sigma_3^-}\,\theta^2+27\,{\sigma_3^-}\,\theta+27
            \,{\sigma_3^-}\right)+3\left({\sigma_3^-}-18\right)\,\theta^2
          -9\left({\sigma_3^-}+6\right)\,\theta}\over{3\left( 4\,\theta^2+12
            \,\theta+9\right) }}
      \frac{\Delta a}{a}
      \\
      & - &
      {{-48\left({\sigma_3^-}+5\right)\,\theta^5+\sqrt{\theta
          }\,\left(144\,{\sigma_3^-}\,\theta^4+780\,{\sigma_3^-}\,
            \theta^3+2214\,{\sigma_3^-}\,\theta^2+1944\,{\sigma_3^-}\,
            \theta+1215\,{\sigma_3^-}\right)}\over{48\left( 2\,\theta
            ^4+12\,\theta^3+27\,\theta^2+27\,\theta+12\right) }}
      \left( \frac{\Delta a}{a} \right)^2
      \\
      & -&
      {{- 4\left(91\,{\sigma_3^-}+342
          \right)\,\theta^4-54\left(11\,{\sigma_3^-}+56\right)\,
          \theta^3-216\left(2\,{\sigma_3^-}+15\right)\,\theta^2-81
          \left(5\,{\sigma_3^-}+18\right)\,\theta}\over{48\left( 2\,\theta
            ^4+12\,\theta^3+27\,\theta^2+27\,\theta+12\right)}}
      \left( \frac{\Delta a}{a} \right)^2
    \end{eqnarray*}
    {\normalsize Above are the values of $\theta^+$ in the
      cases~\textbf{(S)}, \textbf{(P)}, \textbf{(L)}
      and~\textbf{(p)}.}
  \end{landscape}
}

{\small{

    \bibliographystyle{abbrv}

    \bibliography{Eulero}

  }}
\end{document}